\documentclass[10pt]{article}
\usepackage{amsmath}
\usepackage{amsfonts}
\usepackage{amssymb}
\usepackage{graphicx}
\usepackage{mathrsfs}
\usepackage{xcolor}
\usepackage{verbatim}
\usepackage{mathrsfs}
\usepackage[body={15.5cm,21cm}, top=3cm]{geometry}
\usepackage{paralist}

\DeclareMathOperator*{\esssup}{ess\,sup}
\allowdisplaybreaks[4]
\providecommand{\U}[1]{\protect \rule{.1in}{.1in}}
\newtheorem{theorem}{Theorem}[section]

\newenvironment{proof}[1][Proof]{\noindent \textbf{#1.} }{\  \rule{0.5em}{0.5em}}
\allowdisplaybreaks[2]
\begin{document}

\title{Existence of global solutions for multi-dimensional coupled FBSDEs with diagonally quadratic generators}
\author{Yuyang Chen\thanks{School of Mathematical Sciences, Shanghai Jiao Tong University, China (cyy0032@sjtu.edu.cn)}
\and
Peng Luo \thanks{School of Mathematical Sciences, Shanghai Jiao Tong University, China (peng.luo@sjtu.edu.cn). Financial support from the National Natural Science Foundation of China (Grant No. 12101400) is gratefully acknowledged.}}

\maketitle
\begin{abstract}
The present paper is devoted to study multi-dimensional coupled FBSDEs with diagonally quadratic generators. Relying on a comparison result obtained in \cite{19}, we provide conditions under which there exists a global solution. As a byproduct, we further give a comparison result for this global solution.
\end{abstract}


\textbf{MSC-classification}: 60H10, 60H30.

\section{Introduction}
In this work, we consider the following system of forward and backward stochastic differential equations
\begin{equation}\label{eq:FBSDE}
\begin{cases}
X_{t}&=x+\int_{0}^{t}b(s,X_{s},Y_{s})ds+\int_{0}^{t}\sigma(s,X_{s})dW_{s}\\
Y_{t}&=h(X_{T})+\int_{t}^{T}g(s,X_{s},Y_{s},Z_{s})ds-\int_{t}^{T}Z_{s}dW_{s},\quad t\in[0,T]
\end{cases}
\end{equation}
where $W$ is a multi-dimensional Brownian motion on a probability space, $x$ is the initial condition, $T>0$ is a fixed finite time horizon, and $b,\sigma,h,g$ are functions. We provide coniditions which guarantee the existence of global solutions in the case where $Y$ is multi-dimensional and $g$ is diagonally quadratic in $z$.

Solvability of coupled FBSDEs with Lipschitz generators has been well studied in the literature. Antonelli \cite{Anto} obtained the first solvability result of a coupled FBSDE over a small time horizon, which is further investigated in Pardoux and Tang \cite{PT}. Ma et al. \cite{6} proposed a four step scheme method to solve couplded Markovian FBSDEs (see also \cite{7}). Their method allows to obtain solvability result for arbitrarily large time horizon. Hu and Peng \cite{1} introduces a continuation method, which is further developed in \cite{2,3}. This approach is able to provide global solvability for coupled non-Markovian FBSDEs under a so-called monotonicity conidition. Combining contraction mapping method and a four step scheme method, Delarue \cite{5} obtained the existence and uniqueness of global solutions for coupled FBSDEs with non-degenerate diffusion processes. Antonelli and Hamad\`ene \cite{8} investigated the solvability of coupled FBSDEs with continuous monotone coefficients. Recently, Ma et al. \cite{MWZZ} established a unified approach to study the wellposedness of non-Markovian coupled FBSDE, where the concept of decoupling field has been introduced. This approach is further extended by Fromm and Imkeller \cite{15}.

On the other hand, BSDEs and FBSDEs with quadratic growth naturally arise in mathematical finance, see e.g. Horst et al. \cite{9}, Kramkov and Pulido \cite{10} and Bielagk et al. \cite{11}. One-dimensional quadratic BSDE was first studied by Kobylanski \cite{12}, where the existence, comparison and stability results are obtained. Briand and Hu \cite{BH,BH1} further studied quadratic BSDEs for unbounded terminal conditions. Barrieu and El Karoui \cite{13} considered quadratic BSDEs in terms of some general quadratic semimartingales. We further refer to Tevzadz \cite{Te}, Hu and Tang \cite{21}, Jamneshan et al. \cite{JKL}, Luo \cite{Luo}, Cheridito and Nam \cite{CN} and Xing and \v{Z}itkovi\'{c} \cite{XZ} for recent developments in solvability of multi-dimensional quadratic BSDEs.

The main scope of this paper is to consider the solvability of system \eqref{eq:FBSDE} where $Y$ is multi-dimensional and $g$ has diagonally quadratic growth in $z$. When $Y$ is one dimensional and $g$ is quadratic in $z$, Antonelli and Hamad\`{e}ne \cite{8} obtained the existence of a global solution for system \eqref{eq:FBSDE}. Using the technique of decoupling field and in a Markovian case, Fromm and Imkeller \cite{15} obtained the existence and uniqueness of a local solution and provided an extension to maximally solvable time horizon for generators which are locally Lipschitz in $z$. When $\sigma$ is independent of $X$, Luo and Tangpi \cite{14} obtained the existence and uniqueness of a local solution for generators which can be separated into a quadratic and subquadratic part, while Kupper et al. \cite{17} obtained global solvability for generators with arbitrary growth in $z$ in a Markovian setting. More recently, Jackson\cite{16} established global solvability for a type of  quadratic Markovian FBSDEs. Compared with all these results, we study global solution for system \eqref{eq:FBSDE} in a multi-dimensional and non-Markovian framework. Relying on a comparison result recently established in Luo \cite{19}, we construct a sequence of processes and obtain some delicated a prior estimates. By making ample use of properties of BMO martingales, we further show the convergence of this sequence of processes in suitable spaces, which yields the existence of a global solution of system \eqref{eq:FBSDE}. Finally, we present a comparison result for this global solution. The paper is organized as follows. In the next section, we state our setups and main results. Theorem 2.1 provides the existence of global solutions of multi-dimensional FBSDEs with diagonally quadratic generators. A comparison result for this global solution is given in Theorem 2.2. 
\section{Main results}
Consider the following forward and backward stochastic differential equations(FBSDEs):
\begin{equation}\label{eq:FBSDE_main}
\begin{cases}
X_{t}^{i}&=x_{0}^{i}+\int_{0}^{t}b^{i}(s,X_{s},Y_{s})ds+\int_{0}^{t}\sigma^{i}(s,X_{s})dW_{s},\\
Y_{t}^{i}&=h^{i}(X_{T})+\int_{t}^{T}g^{i}(s,X_{s},Y_{s},Z_{s}^{i})ds-\int_{t}^{T}Z_{s}^{i}dW_{s},
\end{cases}
\end{equation}
for $i=1,\cdots,n,~t\in[0,T]$. Here, $W=(W_{t})_{t\geq0}$ is a $d$-dimensional standard Brownian motion defined on some probability space $(\Omega,\mathcal{F},P)$. Denote by $\{\mathcal{F}_{t},0\leq t\leq T\}$ the augmented natural filtration of $W$. Equalities and inequalities between random variables and processes are understood in the $P$-a.s. and $P\otimes dt$-a.e. sense, respectively. The Euclidean norm is denoted by $|\cdot|$ and for a random variable $\eta$, $\|\eta\|_{\infty}$ denotes the $L^{\infty}$-norm of $\eta$, i.e., $\|\eta\|_{\infty}:=\mathop{\esssup}\limits_{\omega}|\eta(\omega)|$. For $x,y\in\mathbb{R}^{n},x\leq y$ is understood component-wisely, i.e., $x\leq y$ if and only if $x^{i}\leq y^{i}$ for all $i=1,\cdots,n$. For $p\geq2$, we denote by
\begin{itemize}
  \item $\mathcal{H}^{p}(\mathbb{R}^{n})$\ the space of $n$-dimensional predictable processes $Y$ on [0,T] such that\\
  \centerline{$\|Y\|_{\mathcal{H}^{p}(\mathbb{R}^{n})}:=E\left[\left(\int_{0}^{T}|Y_{t}|^{2}dt\right)^{\frac{p}{2}}\right]^{\frac{1}{p}}<\infty$;}\\
  \item $\mathcal{S}^{p}(\mathbb{R}^{n})$\ the space of $n$-dimensional predictable processes $Y$ on [0,T] such that\\
  \centerline{$\|Y\|_{\mathcal{S}^{p}(\mathbb{R}^{n})}:=E\left[\left(\sup\limits_{0\leq t\leq T}|Y_{t}|^{2}\right)^{\frac{p}{2}}\right]^{\frac{1}{p}}<\infty$;}\\
  \item $\mathcal{S}^{\infty}(\mathbb{R}^{n})$\ the space of $n$-dimensional predictable processes $Y$ on [0,T] such that\\
  \centerline{$\|Y\|_{\mathcal{S}^{\infty}(\mathbb{R}^{n})}:=\left\|\sup\limits_{0\leq t\leq T}|Y_{t}|\right\|_{\infty}<\infty$.}
\end{itemize}
Let $\mathcal{T}$ be the set of all stopping times with values in $[0,T]$. For any uniformly integrable martingale $M$ with $M_{0}=0$, we set
\begin{align*}
\|M\|_{BMO_{2}}:=\sup_{\tau\in\mathcal{T}}\|E[|M_{T}-M_{\mathcal{\tau}}|^{2}|\mathcal{F}_{\mathcal{T}}]^{\frac{1}{2}}\|_{\infty}.
\end{align*}
The class $\{M:\|M\|_{BMO_{2}}<\infty\}$ is denoted by $BMO$. For $(\alpha\cdot W)_{t}:=\int_{0}^{t}\alpha_{s}dW_{s}$ in $BMO$, the corresponding stochastic exponential is denoted by $\mathcal{E}_{t}(\alpha\cdot W)$.\\
In this paper, we make the following assumptions. Let $C$ be a positive constant.
\begin{flushleft}
$(\mathscr{A}1)$The function $b:=(b^{1},\cdots,b^{n})^{*}:\Omega\times[0,T]\times\mathbb{R}^{n}\times\mathbb{R}^{n}\rightarrow\mathbb{R}^{n}$ satisfies that $b(\cdot,x,y)$ is adapted for each $x\in\mathbb{R}^{n}$ and $y\in\mathbb{R}^{n}$. It holds that
\begin{align*}
&|b(t,x,y)|\leq C(1+|x|+|y|),\\
&|b(t,x,y)-b(t,\overline{x},\overline{y})|\leq C(|x-\overline{x}|+|y-\overline{y}|)
\end{align*}
for $x,\overline{x},y,\overline{y}\in\mathbb{R}^{n}$.
\end{flushleft}
\begin{flushleft}
$(\mathscr{A}2)$The function $\sigma:=(\sigma^{1},\cdots,\sigma^{n})^{*}:\Omega\times[0,T]\times\mathbb{R}^{n}\rightarrow\mathbb{R}^{n\times d}$ satisfies that $\sigma(\cdot,x)$ is adapted for each $x\in\mathbb{R}^{n}$. It holds that
\begin{align*}
&|\sigma(t,x)|\leq C(1+|x|),\\
&|\sigma(t,x)-\sigma(t,\overline{x})|\leq C|x-\overline{x}|
\end{align*}
for $x,\overline{x}\in\mathbb{R}^{n}$.
\end{flushleft}
\begin{flushleft}
$(\mathscr{A}3)$The function $h:=(h^{1},\cdots,h^{n})^{*}:\Omega\times\mathbb{R}^{n}\rightarrow\mathbb{R}^{n}$ satisfies that $h$ is $\mathcal{F}_{T}$-measurable. It holds that
\begin{align*}
&|h(x)|\leq C,\\
&|h(x)-h(\overline{x})|\leq C|x-\overline{x}|
\end{align*}
for $x,\overline{x}\in\mathbb{R}^{n}$.
\end{flushleft}
\begin{flushleft}
$(\mathscr{A}4)$For $i=1,\cdots,n$, the function $g^{i}:\Omega\times[0,T]\times\mathbb{R}^{n}\times\mathbb{R}^{n}\times\mathbb{R}^{d}\rightarrow\mathbb{R}$ satisfies that $g^{i}(\cdot,x,y,z^{i})$ is adapted for each $x\in\mathbb{R}^{n}$, $y\in\mathbb{R}^{n}$ and $z^{i}\in\mathbb{R}^{d}$. It holds that
\begin{align*}
&|g^{i}(t,x,y,z^{i})|\leq C(1+|y|+|z^{i}|^{2}),\\
&|g^{i}(t,x,y,z^{i})-g^{i}(t,\overline{x},\overline{y},\overline{z}^{i})|\leq C|x-\overline{x}|+C|y-\overline{y}|+C(1+|z^{i}|+|\overline{z}^{i}|)|z^{i}-\overline{z}^{i}|
\end{align*}
for $x,\overline{x},y,\overline{y}\in\mathbb{R}^{n}$, $z^{i},\overline{z}^{i}\in\mathbb{R}^{d}$.
\end{flushleft}
\begin{flushleft}
$(\mathscr{A}5)$For $t\in[0,T],~i=1,\cdots,n$, it holds that $b^{i}(t,x,y)\leq b^{i}(t,\overline{x},\overline{y})$ for any $x,\overline{x},y,\overline{y}$ satisfying $x^{i}=\overline{x}^{i},x^{j}\leq \overline{x}^{j},j\neq i,y\leq \overline{y}$.
\end{flushleft}
\begin{flushleft}
$(\mathscr{A}6)$For $t\in[0,T]$, it holds that $h(x)\leq h(\overline{x})$ for any $x,\overline{x}$ satisfying $x\leq \overline{x}$.
\end{flushleft}
\begin{flushleft}
$(\mathscr{A}7)$For $t\in[0,T],~i=1,\cdots,n$, it holds that $g^{i}(t,x,y,z^{i})\leq g^{i}(t,\overline{x},\overline{y},z^{i})$ for any $x,\overline{x},y,\overline{y}$ satisfying $y^{i}=\overline{y}^{i},y^{j}\leq \overline{y}^{j},j\neq i,x\leq \overline{x}$.
\end{flushleft}

\begin{theorem}
Let assumptions $(\mathscr{A}1)-(\mathscr{A}7)$ be satisfied, then there exists a solution $(X,Y,Z)$ of \eqref{eq:FBSDE_main} such that $(X,Y,Z\cdot W)\in \mathcal{S}^{p}(\mathbb{R}^{n})\times \mathcal{S}^{\infty}(\mathbb{R}^{n})\times BMO$ for any $p\geq2$. Besides this solution is the minimal one, in the sense that if $(X',Y',Z')$ is another solution of \eqref{eq:FBSDE_main}, for any $t\in[0,T]$, we have
\begin{align*}
X_{t}\leq X'_{t},~Y_{t}\leq Y'_{t}.
\end{align*}
\end{theorem}
\begin{proof}
We divide the proof into three steps:\\
$\mathbf{Step\ 1}$: For any $t\in[0,T],~i=1,\cdots,n$, we consider the following BSDE and SDE:
\begin{align}
&U_{t}^{i}=C+\int_{t}^{T}C(1+|U_{s}|+|V_{s}^{i}|^{2})ds-\int_{t}^{T}V_{s}^{i}dW_{s},\\
&S_{t}^{i}=x_{0}^{i}+\int_{0}^{t}C(1+|S_{s}|+|U_{s}|)ds+\int_{0}^{t}\sigma^{i}(s,S_{s})dW_{s}.
\end{align}
There exists a unique solution $(U,V)$ in (2) such that $(U,V)\in\mathcal{S}^{\infty}(\mathbb{R}^{n})\times BMO$ from \cite[Theorem 2.3]{21}. Moreover,
\begin{align*}
\|U\|_{\mathcal{S}^{\infty}(\mathbb{R}^{n})}+\|V\cdot W\|_{BMO}\leq K
\end{align*}
where $K$ is a positive constant only depending on $C$ and $T$. Also, the SDE (3) has a solution $S\in\mathcal{S}^{p}(\mathbb{R}^{n})$ due to the SDE theory.\\
Then we consider $(Y^{(0)},Z^{(0)})$ is the solution of
\begin{align}
Y_{t}^{i,(0)}=-C-\int_{t}^{T}C(1+|Y_{s}^{(0)}|+|Z_{s}^{i,(0)}|^{2})ds-\int_{t}^{T}Z_{s}^{i,(0)}dW_{s}.
\end{align}
Noting the assumption $(\mathscr{A}4)$, we can obtain $Y_{t}^{(0)}\leq U_{t}^{(0)}$ for any $t\in[0,T]$ from \cite[Theorem 2.2]{19}. Without loss of generality, assume $\|Y^{(0)}\|_{\mathcal{S}^{\infty}(\mathbb{R}^{n})}+\|Z^{(0)}\cdot W\|_{BMO}\leq K$.\\
Similarly, we consider $X^{(0)}$ is the solution of
\begin{align}
X_{t}^{i,(0)}=x_{0}^{i}+\int_{0}^{t}b^{i}(s,X_{s}^{(0)},Y_{s}^{(0)})ds+\int_{0}^{t}\sigma^{i}(s,X_{s}^{(0)})dW_{s}.
\end{align}
Noting the assumption $(\mathscr{A}1)$, we can obtain $X_{t}^{(0)}\leq S_{t}^{(0)}$ for any $t\in[0,T]$ from \cite[Theorem 1.1]{18}.\\
Now we construct some FBSDEs. For any $t\in[0,T],~k\geq1,~i=1,\cdots,n$,
\begin{align}
&Y_{t}^{i,(k)}=h^{i}(X_{T}^{(k-1)})+\int_{t}^{T}g^{i}(s,X_{s}^{(k-1)},Y_{s}^{(k)},Z_{s}^{i,(k)})ds-\int_{t}^{T}Z_{s}^{i,(k)}dW_{s},\\
&X_{t}^{i,(k)}=x_{0}^{i}+\int_{0}^{t}b^{i}(s,X_{s}^{(k)},Y_{s}^{(k)})ds+\int_{0}^{t}\sigma^{i}(s,X_{s}^{(k)})dW_{s}.
\end{align}
If it holds for any $t\in[0,T]$ that
\begin{align*}
X_{t}^{(k-1)}\leq X_{t}^{(k)}\leq S_{t},
\end{align*}
we can immediately get
\begin{align*}
&h(X_{T}^{(k-1)})\leq h(X_{T}^{(k)})\leq C,\\
&g^{i}(t,X_{t}^{(k-1)},y,z^{i})\leq g^{i}(t,X_{t}^{(k)},y,z^{i})\leq C(1+|y|+|z^{i}|^{2})
\end{align*}
due to the assumptions $(\mathscr{A}6)$$(\mathscr{A}7)$.
It follows from \cite[Theorem 2.2]{19} that
\begin{align*}
Y_{t}^{(k)}\leq Y_{t}^{(k+1)}\leq U_{t},~t\in[0,T].
\end{align*}
Similarly, when
\begin{align*}
Y_{t}^{(k)}\leq Y_{t}^{(k+1)}\leq U_{t}
\end{align*}
holds for any $t\in[0,T]$, we can obtain
\begin{align*}
b(t,x,Y_{t}^{(k)})\leq b(t,x,Y_{t}^{(k+1)})\leq C(1+|x|+|U_{t}|)
\end{align*}
due to the assumption $(\mathscr{A}5)$. Using \cite[Theorem 1.1]{18}, we will get
\begin{align*}
X_{t}^{(k)}\leq X_{t}^{(k+1)}\leq S_{t},~t\in[0,T].
\end{align*}
Since it is easy to check that
\begin{align*}
X_{t}^{(0)}\leq X_{t}^{(1)}\leq S_{t},~Y_{t}^{(0)}\leq Y_{t}^{(1)}\leq U_{t}
\end{align*}
are satisfied for any $t\in[0,T]$, then for any $k\geq0$, we will have
\begin{align*}
X_{t}^{(k)}\leq X_{t}^{(k+1)}\leq S_{t},~Y_{t}^{(k)}\leq Y_{t}^{(k+1)}\leq U_{t},~t\in[0,T]
\end{align*}
by induction. By the monotonic convergence theorem, we note that there exists two adapted processes $X$ and $Y$ such that $\|X^{(k)}-X\|_{\mathcal{H}^{p}(\mathbb{R}^{n})}\rightarrow0,\|Y^{(k)}-Y\|_{\mathcal{H}^{p}(\mathbb{R}^{n})}\rightarrow0$ as $k\rightarrow\infty$ for any $p\geq2$. What's more, 
\begin{align*}
E\left[\int_{0}^{T}\left|X_{t}^{(k)}-X_{t}\right|^{p}dt\right]\rightarrow0,~E\left[\int_{0}^{T}\left|Y_{t}^{(k)}-Y_{t}\right|^{p}dt\right]\rightarrow0
\end{align*}
as $k\rightarrow\infty$ for any $p\geq2$.\\
\\
$\mathbf{Step\ 2}$: We prove that $(X^{(k)},Y^{(k)},Z^{(k)})$ converges to $(X,Y,Z)$ in $\mathcal{S}^{p}(\mathbb{R}^{n})\times\mathcal{S}^{p}(\mathbb{R}^{n})\times\mathcal{H}^{p}(\mathbb{R}^{n\times d})$ as $k\rightarrow\infty$ for any $p\geq2$.\\
Since for any $p\geq2,~k,h\geq0$ it holds that
\begin{align*}
E&\left[\left(\sup\limits_{0\leq t\leq T}|\int_{0}^{t}(b(s,X_{s}^{(k)},Y_{s}^{(k)})-b(s,X_{s}^{(h)},Y_{s}^{(h)}))ds|^{2}\right)^{\frac{p}{2}}\right]^{\frac{1}{p}}\\
\leq&E\left[\left(\sup\limits_{0\leq t\leq T}\int_{0}^{t}|b(s,X_{s}^{(k)},Y_{s}^{(k)})-b(s,X_{s}^{(h)},Y_{s}^{(h)})|^{2}ds\right)^{\frac{p}{2}}\right]^{\frac{1}{p}}\\
\leq&E\left[\left(\int_{0}^{T}(C|X_{s}^{(k)}-X_{s}^{(h)}|+C|Y_{s}^{(k)}-Y_{s}^{(h)}|)^{2}ds\right)^{\frac{p}{2}}\right]^{\frac{1}{p}}\\
\leq&CE\left[\left(\int_{0}^{T}|X_{s}^{(k)}-X_{s}^{(h)}|^{2}ds\right)^{\frac{p}{2}}\right]^{\frac{1}{p}}+CE\left[\left(\int_{0}^{T}|Y_{s}^{(k)}-Y_{s}^{(h)}|^{2}ds\right)^{\frac{p}{2}}\right]^{\frac{1}{p}}\\
\end{align*}
and
\begin{align*}
E\left[\left(\sup\limits_{0\leq t\leq T}|\int_{0}^{t}(\sigma(s,X_{s}^{(k)})-\sigma(s,X_{s}^{(h)})dW_{s}|^{2}\right)^{\frac{p}{2}}\right]^{\frac{1}{p}}
\leq&C_{p}^{\frac{1}{p}}E\left[\left(\int_{0}^{T}|\sigma(s,X_{s}^{(k)})-\sigma(s,X_{s}^{(h)})|^{2}ds\right)^{\frac{p}{2}}\right]^{\frac{1}{p}}\\
\leq&C_{p}^{\frac{1}{p}}CE\left[\left(\int_{0}^{T}|X_{s}^{(k)}-X_{s}^{(h)}|^{2}ds\right)^{\frac{p}{2}}\right]^{\frac{1}{p}}\\
\end{align*}
where $C_{p}$ is the coefficient in the Burkholder-Davis-Gundy inequality, which implies $\|X^{(k)}-X\|_{\mathcal{S}^{p}(\mathbb{R}^{n})}\rightarrow0$ as $k\rightarrow\infty$ for any $ p\geq2$. Let $k\rightarrow\infty$ in (8), we have
\begin{align*}
X_{t}^{i}=x_{0}^{i}+\int_{0}^{t}b^{i}(s,X_{s},Y_{s})ds+\int_{0}^{t}\sigma^{i}(s,X_{s})dW_{s}.
\end{align*}
Then we discuss the convergence of $Y^{(k)}$ and $Z^{(k)}$. Because of the assumption ($\mathscr{A}$4), we can denote
\begin{align*}
&g^{i}(t,X_{t}^{(k-1)},Y_{t}^{(k)},Z_{t}^{i,(k)})-g^{i}(t,X_{t}^{(h-1)},Y_{t}^{(h)},Z_{t}^{i,(h)})\\
=&\alpha_{t}^{i}(k,h)(X_{t}^{(k-1)}-X_{t}^{(h-1)})+\beta_{t}^{i}(k,h)(Y_{t}^{(k)}-Y_{t}^{(h)})+
\gamma_{t}^{i}(k,h)(Z_{t}^{i,(k)}-Z_{t}^{i,(h)})
\end{align*}
and
\begin{align*}
|\alpha_{t}^{i}(k,h)|\leq C,~|\beta_{t}^{i}(k,h)|\leq C,~|\gamma_{t}^{i}(k,h)|\leq C(1+|Z_{t}^{i,(k)}|+|Z_{t}^{i,(h)}|)
\end{align*}
for any $k,h\geq1,~i=1,\cdots,n$. We note that $W_{t}(i,k,h):=W_{t}-\int_{0}^{t}\gamma_{s}^{i}(k,h)ds$ is a Brownian motion under the equivalent probability measure $P^{i,k,h}$ defined by
\begin{align*}
dP^{i,k,h}:=\mathcal{E}(\gamma^{i}(k,h)\cdot W)_{0}^{T}dP.
\end{align*}
Since $\int_{0}^{t}\gamma_{s}^{i}(k,h)dW_{s}$ is a $BMO$ martingale, we can deduce that there exists $p_{1},q_{1}>1$ depending on $\gamma^{i}(k,h)$, such that
\begin{align*}
E_{\tau}\left(\frac{dP^{i,k,h}}{dP}\right)^{p_{1}}<\infty,~E_{\tau}^{i,k,h}\left(\frac{dP}{dP^{i,k,h}}\right)^{p_{2}}<\infty
\end{align*}
from \cite[Theorem 3.1]{20}, where $E_{\tau}$ is a conditional expectation for an arbitrary stopping time $\tau$. We let $q_{1},q_{2}>1$ satisfying
\begin{align*}
\frac{1}{p_{1}}+\frac{1}{q_{1}}=1,~\frac{1}{p_{2}}+\frac{1}{q_{2}}=1.
\end{align*}
Denote the conditional expectation with respect to $P$ by $E_{t}$ and $P^{i,k,h}$ by $E_{t}^{i,k,h}$ and use It\^{o}'s formula, for any $p\geq2,~i=1,\cdots,n$, we have
\begin{align*}
&|Y_{t}^{i,(k)}-Y_{t}^{i,(h)}|^{p}+\frac{1}{2}p(p-1) E_{t}^{i,k,h}\int_{t}^{T}|Y_{s}^{i,(k)}-Y_{s}^{i,(h)}|^{p-2}|Z_{s}^{i,(k)}-Z_{s}^{i,(h)}|^{2}ds\\
&=E_{t}^{i,k,h}|h^{i}(X_{T}^{(k-1)})-h^{i}(X_{T}^{(h-1)})|^{p}\\
&+E_{t}^{i,k,h}\int_{t}^{T}p|Y_{s}^{i,(k)}-Y_{s}^{i,(h)}|^{p-2}(Y_{s}^{i,(k)}-Y_{s}^{i,(h)})[\alpha_{s}^{i}(k,h)(X_{s}^{(k-1)}-X_{s}^{(h-1)})+\beta_{s}^{i}(k,h)(Y_{s}^{(k)}-Y_{s}^{(h)})]ds.
\end{align*}
First, we calculate
\begin{align*}
&E_{t}^{i,k,h}\int_{t}^{T}p |Y_{s}^{i,(k)}-Y_{s}^{i,(h)}|^{p-2}(Y_{s}^{i,(k)}-Y_{s}^{i,(h)})\alpha_{s}^{i}(k,h)(X_{s}^{(k-1)}-X_{s}^{(h-1)})ds\\
&\leq pCE_{t}^{i,k,h}\int_{t}^{T}|Y_{s}^{(k)}-Y_{s}^{(h)}|^{p-1}|X_{s}^{(k-1)}-X_{s}^{(h-1)}|ds\\
&\leq (p-1)CE_{t}^{i,k,h}\int_{t}^{T}|Y_{s}^{(k)}-Y_{s}^{(h)}|^{p}ds+CE_{t}^{i,k,h}\int_{t}^{T}|X_{s}^{(k-1)}-X_{s}^{(h-1)}|^{p}ds
\end{align*}
and
\begin{align*}
&E_{t}^{i,k,h}\int_{t}^{T}p |Y_{s}^{i,(k)}-Y_{s}^{i,(h)}|^{p-2}(Y_{s}^{i,(k)}-Y_{s}^{i,(h)})\beta_{s}^{i}(k,h)(Y_{s}^{(k)}-Y_{s}^{(h)})ds\\
&\leq pCE_{t}^{i,k,h}\int_{t}^{T}|Y_{s}^{(k)}-Y_{s}^{(h)}|^{p}ds
\end{align*}
for any $p\geq2,~k,h\geq1,~i=1,\cdots,n$, which implies
\begin{align*}
|Y_{t}^{i,(k)}-Y_{t}^{i,(h)}|^{p}\leq &C^{p}E_{t}^{i,k,h}|X_{T}^{(k-1)}-X_{T}^{(h-1)}|^{p}+(2p-1)CE_{t}^{i,k,h}\int_{t}^{T}|Y_{s}^{(k)}-Y_{s}^{(h)}|^{p}ds\\
&+CE_{t}^{i,k,h}\int_{t}^{T}|X_{s}^{(k-1)}-X_{s}^{(h-1)}|^{p}ds.
\end{align*}
Then we show that for any $p\geq2,~k,h\geq1,~i=1,\cdots,n$,
\begin{align*}
E^{i,k,h}\left[\sup_{t\in[0,T]}E_{t}^{i,k,h}\int_{t}^{T}|Y_{s}^{(k)}-Y_{s}^{(h)}|^{p}ds\right]&\leq2E^{i,k,h}\left[\left(\int_{0}^{T}|Y_{s}^{(k)}-Y_{s}^{(h)}|^{p}ds\right)^{2}\right]^{\frac{1}{2}}\\
&=2E\left[\left(\frac{dP^{i,k,h}}{dP}\right)\left(\int_{0}^{T}|Y_{s}^{(k)}-Y_{s}^{(h)}|^{p}ds\right)^{2}\right]^{\frac{1}{2}}\\
&\leq2E\left[\left(\frac{dP^{i,k,h}}{dP}\right)^{p_{1}}\right]^{\frac{1}{2p_{1}}}E\left[\left(\int_{0}^{T}|Y_{s}^{(k)}-Y_{s}^{(h)}|^{p}ds\right)^{2q_{1}}\right]^{\frac{1}{2q_{1}}}\\
&\leq2E\left[\left(\frac{dP^{i,k,h}}{dP}\right)^{p_{1}}\right]^{\frac{1}{2p_{1}}}E\left[\left(\int_{0}^{T}|Y_{s}^{(k)}-Y_{s}^{(h)}|^{2pq_{1}}ds\right)\right]^{\frac{1}{2q_{1}}}
\end{align*}
holds due to Doob's $L^{p}$-inequality and H\"older's inequality.\\
Similarly, we have
\begin{align*}
&E^{i,k,h}\left[\sup_{t\in[0,T]}E_{t}^{i,k,h}\int_{t}^{T}|X_{s}^{(k-1)}-X_{s}^{(h-1)}|^{p}ds\right]\leq2E\left[\left(\frac{dP^{i,k,h}}{dP}\right)^{p_{1}}\right]^{\frac{1}{2p_{1}}}E\left[\left(\int_{0}^{T}|X_{s}^{(k-1)}-X_{s}^{(h-1)}|^{2pq_{1}}ds\right)\right]^{\frac{1}{2q_{1}}},\\
&E^{i,k,h}\left[\sup_{t\in[0,T]}E_{t}^{i,k,h}\left|X_{T}^{(k-1)}-X_{T}^{(h-1)}\right|^{p}\right]\leq 2E\left[\left(\frac{dP^{i,k,h}}{dP}\right)^{p_{1}}\right]^{\frac{1}{2p_{1}}}E\left[\left|X_{T}^{(k-1)}-X_{T}^{(h-1)}\right|^{2pq_{1}}\right]^{\frac{1}{2q_{1}}}.
\end{align*}
So it holds that for any $p\geq2,~k,h\geq1,~i=1,\cdots,n$,
\begin{align*}
E^{i,k,h}\left[\sup_{t\in[0,T]}\left|Y_{t}^{i,(k)}-Y_{t}^{i,(h)}\right|^{p}\right]
&\leq 2C^{p}E\left[\left(\frac{dP^{i,k,h}}{dP}\right)^{p_{1}}\right]^{\frac{1}{2p_{1}}}E\left[\left|X_{T}^{(k-1)}-X_{T}^{(h-1)}\right|^{2pq_{1}}\right]^{\frac{1}{2pq_{1}}}\\
&+2(2p-1)CE\left[\left(\frac{dP^{i,k,h}}{dP}\right)^{p_{1}}\right]^{\frac{1}{2p_{1}}}E\left[\left(\int_{0}^{T}|Y_{s}^{(k)}-Y_{s}^{(h)}|^{2pq_{1}}ds\right)\right]^{\frac{1}{2q_{1}}}\\
&+2CE\left[\left(\frac{dP^{i,k,h}}{dP}\right)^{p_{1}}\right]^{\frac{1}{2p_{1}}}E\left[\left(\int_{0}^{T}|X_{s}^{(k-1)}-X_{s}^{(h-1)}|^{2pq_{1}}ds\right)\right]^{\frac{1}{2q_{1}}}.
\end{align*}
Since we have that
\begin{align*}
E\left[\sup_{t\in[0,T]}\left|Y_{t}^{i,(k)}-Y_{t}^{i,(h)}\right|^{p}\right]
&=E^{i,k,h}\left[\left(\frac{dP}{dP^{i,k,h}}\right)\sup_{t\in[0,T]}\left|Y_{t}^{i,(k)}-Y_{t}^{i,(h)}\right|^{p}\right]\\
&\leq E^{i,k,h}\left[\left(\frac{dP}{dP^{i,k,h}}\right)^{p_{2}}\right]^{\frac{1}{p_{2}}}E^{i,k,h}\left[\sup_{t\in[0,T]}\left|Y_{t}^{i,(k)}-Y_{t}^{i,(h)}\right|^{pq_{2}}\right]^{\frac{1}{q_{2}}}
\end{align*}
and
\begin{align*}
&E\left(\frac{dP^{i,k,h}}{dP}\right)^{p_{1}}<\infty,~E^{i,k,h}\left(\frac{dP}{dP^{i,k,h}}\right)^{p_{2}}<\infty,\\
&\lim_{k,h\rightarrow\infty}E\left[\int_{0}^{T}\left|X_{t}^{(k)}-X_{t}^{(h)}\right|^{p}dt\right]=0,~\lim_{k,h\rightarrow\infty}E\left[\int_{0}^{T}\left|Y_{t}^{(k)}-Y_{t}^{(h)}\right|^{p}dt\right]=0
\end{align*}
for any $p\geq2$, we can obtain
\begin{align*}
\lim\limits_{k,h\rightarrow\infty}E\left[\sup\limits_{t\in[0,T]}\left|Y_{t}^{(k)}-Y_{t}^{(h)}\right|^{p}\right]=0
\end{align*}
for any $p\geq2$. So $\|Y^{(k)}-Y\|_{\mathcal{S}^{p}(\mathbb{R}^{n})}\rightarrow0$ as $k\rightarrow\infty$ for any $p\geq2$.\\
We show it holds that $(Y^{(k)},Z^{(k)}\cdot W)\in\mathcal{S}^{\infty}(\mathbb{R}^{n})\times BMO$ for any $k\geq0$.\\
Since we have
\begin{align*}
Y_{t}^{(0)}\leq Y_{t}^{(k)}\leq U_{t},~t\in[0,T]
\end{align*}
and $Y^{(k)}$ is continuous in $\mathcal{S}^{p}(\mathbb{R}^{n})$ for any $p\geq2,~k\geq0$, it holds that
\begin{align}
\|Y^{(k)}\|_{\mathcal{S}^{\infty}(\mathbb{R}^{n})}\leq \max\{\|Y^{(0)}\|_{\mathcal{S}^{\infty}(\mathbb{R}^{n})},\|U\|_{\mathcal{S}^{\infty}(\mathbb{R}^{n})}\}\leq K,
\end{align}
which implies $Y^{(k)}$ is in $S^{\infty}(\mathbb{R}^{n})$.\\
Besides, define
\begin{align}
\phi(x)=\frac{e^{2C|x|}-2C|x|-1}{4C^{2}}\geq0
\end{align}
and we can calculate that
\begin{align*}
\phi'(x)=\frac{e^{2C|x|}-1}{2C}sgn(x),~\phi''(x)=e^{2C|x|},~\phi''(x)-2C|\phi'(x)|=1.
\end{align*}
By It\^o's formula, it holds that for any $k\geq1,~i=1,\cdots,n$,
\begin{align*}
\phi(Y_{t}^{i,(k)})+\frac{1}{2}E_{t}\int_{t}^{T}\phi''(Y_{s}^{i,(k)})|Z_{s}^{i,(k)}|^{2}ds
&=E_{t}\phi(Y_{T}^{i,(k)})+E_{t}\int_{t}^{T}\phi'(Y_{s}^{i,(k)})g^{i}(X_{s},Y_{s}^{(k)},Z_{s}^{i,(k)})ds\\
&\leq\phi(\|Y_{T}^{i,(k)}\|_{\infty})+E_{t}\int_{t}^{T}\phi'(Y_{s}^{i,(k)})C(1+|Y_{s}^{(k)}|+|Z_{s}^{i,(k)}|^{2})ds,
\end{align*}
and then
\begin{equation}
\begin{aligned}
\phi(Y_{t}^{i,(k)})+\frac{1}{2}E_{t}\int_{t}^{T}|Z_{s}^{i,(k)}|^{2}ds
&\leq\phi(\|Y_{T}^{i,(k)}\|_{\infty})+E_{t}\int_{t}^{T}\phi'(Y_{s}^{i,(k)})C(1+|Y_{s}^{(k)}|)ds\\
&\leq\phi(\|Y_{T}^{(k)}\|_{\infty})+C(T-t)\phi'(\|Y^{(k)}\|_{\mathcal{S}^{\infty}(\mathbb{R}^{n})})(1+\|Y\|_{\mathcal{S}^{\infty}(\mathbb{R}^{n})})\\
&\leq\phi(K)+CT\phi'(K)(1+K),
\end{aligned}
\end{equation}
which implies
\begin{align*}
\|Z^{(k)}\cdot W\|_{BMO}<\infty.
\end{align*}
Then use It\^{o}'s formula for $|Y_{t}^{i,(k)}-Y_{t}^{i,(h)}|^{2}$, we have
\begin{align*}
&|Y_{t}^{i,(k)}-Y_{t}^{i,(h)}|^{2}+\int_{t}^{T}|Z_{s}^{i,(k)}-Z_{s}^{i,(h)}|^{2}ds\\
&=|h^{i}(X_{T}^{(k-1)})-h^{i}(X_{T}^{(h-1)})|^{2}-2\int_{t}^{T}(Y_{s}^{i,(k)}-Y_{s}^{i,(h)})(Z_{s}^{i,(k)}-Z_{s}^{i,(h)})dW_{s}\\
&+2\int_{t}^{T}(Y_{s}^{i,(k)}-Y_{s}^{i,(h)})[\alpha_{s}^{i}(k,h)(X_{s}^{(k-1)}-X_{s}^{(h-1)})+\beta_{s}^{i}(k,h)(Y_{s}^{(k)}-Y_{s}^{(h)})+\gamma_{s}^{i}(k,h)(Z_{s}^{i,(k)}-Z_{s}^{i,(h)})]ds.
\end{align*}
We can get
\begin{align*}
\int_{0}^{T}|Z_{t}^{i,(k)}-Z_{t}^{i,(h)}|^{2}dt
\leq &C^{2}|X_{T}^{(k-1)}-X_{T}^{(h-1)}|^{2}+2\int_{0}^{T}|(Y_{t}^{i,(k)}-Y_{t}^{i,(h)})(Z_{t}^{i,(k)}-Z_{t}^{i,(h)})|dW_{t}\\
&+2CK\int_{0}^{T}|X_{s}^{(k-1)}-X_{s}^{(h-1)}|ds+2CT\sup_{t\in[0,T]}|Y_{t}^{(k)}-Y_{t}^{(h)}|^{2}\\
&+6C\sup_{t\in[0,T]}|Y_{t}^{(k)}-Y_{t}^{(h)}|\int_{0}^{T}(1+|Z_{t}^{i,(k)}|^{2}+|Z_{t}^{i,(h)}|^{2})dt.
\end{align*}
With $C_{r}$ inequality and H\"{o}lder's inequality, for any $p\geq2$, it holds that
\begin{align*}
&E[(\int_{0}^{T}|Z_{s}^{i,(k)}-Z_{s}^{i,(h)}|^{2}ds)^{\frac{p}{2}}]\\
\leq& C_{p}\left\{C^{p}E\left|X_{T}^{(k-1)}-X_{T}^{(h-1)}\right|^{p}+2^{\frac{p}{2}}E\left[\left(\int_{0}^{T}|(Y_{s}^{i,(k)}-Y_{s}^{i,(h)})(Z_{s}^{i,(k)}-Z_{s}^{i,(h)})|^{2}ds\right)^{\frac{p}{4}}\right]\right.\\
&+(2CK)^{\frac{p}{2}}E\left[\left(\int_{0}^{T}|X_{s}^{(k-1)}-X_{s}^{(h-1)}|ds\right)^{\frac{p}{2}}\right]+(2CT)^{\frac{p}{2}}E\left[\left(\sup_{t\in[0,T]}\left|Y_{t}^{(k)}-Y_{t}^{(h)}\right|^{2}\right)^{\frac{p}{2}}\right]\\
&\left.+(6C)^{\frac{p}{2}}E\left[\left(\sup_{t\in[0,T]}\left|Y_{t}^{(k)}-Y_{t}^{(h)}\right|\right)^{\frac{p}{2}}\left(\int_{0}^{T}(1+|Z_{s}^{i,(k)}|^{2}+|Z_{s}^{i,(h)}|^{2})ds\right)^{\frac{p}{2}}\right]\right\}\\
\leq& C'_{p}\left\{E\left|X_{T}^{(k-1)}-X_{T}^{(h-1)}\right|^{p}+E\left[\sup_{t\in[0,T]}\left|Y_{t}^{(k)}-Y_{t}^{(h)}\right|^{p}\right]^{\frac{1}{2}}E\left[\left(\int_{0}^{T}(|Z_{s}^{i,(k)}|^{2}+|Z_{s}^{i,(h)})|^{2})ds\right)^{\frac{p}{2}}\right]^{\frac{1}{2}}\right.\\
&+E\left[\left(\int_{0}^{T}|X_{s}^{(k-1)}-X_{s}^{(h-1)}|^{2}ds\right)^{\frac{p}{4}}\right]+E\left[\left(\sup_{t\in[0,T]}\left|Y_{t}^{(k)}-Y_{t}^{(h)}\right|^{2}\right)^{\frac{p}{2}}\right]\\
&\left.+E\left[\sup_{t\in[0,T]}\left|Y_{t}^{(k)}-Y_{t}^{(h)}\right|^{p}\right]^{\frac{1}{2}}E\left[\left(\int_{0}^{T}(1+|Z_{s}^{i,(k)}|^{2}+|Z_{s}^{i,(h)}|^{2})ds\right)^{p}\right]^{\frac{1}{2}}\right\}.
\end{align*}
Since for any $p\geq2,~k,h\geq0$, we have
\begin{align*}
E\left(\int_{0}^{T}\left|Z_{s}^{i,(k)}\right|^{2}ds\right)^{\frac{p}{2}}<\infty
\end{align*}
and
\begin{align*}
\lim\limits_{k,h\rightarrow\infty}E\left[\sup\limits_{t\in[0,T]}\left|X_{t}^{(k)}-X_{t}^{(h)}\right|^{p}\right]=0,~\lim\limits_{k,h\rightarrow\infty}E\left[\sup\limits_{t\in[0,T]}\left|Y_{t}^{(k)}-Y_{t}^{(h)}\right|^{p}\right]=0,
\end{align*}
then $Z^{(k)}$ is convergent in $\mathcal{H}^{p}(\mathbb{R}^{n\times d})$. So there exists $Z\in\mathcal{H}^{p}(\mathbb{R}^{n\times d})$, such that $\|Z^{(k)}-Z\|_{\mathcal{H}^{p}(\mathbb{R}^{n\times d})}\rightarrow0$ as $k\rightarrow\infty$ for any $p\geq2$.\\
Therefore let $k\rightarrow\infty$ in (7), we have
\begin{align*}
Y_{t}^{i}=h(X_{T}^{i})+\int_{0}^{T}g^{i}(s,X_{s},Y_{s},Z_{s}^{i})ds-\int_{0}^{T}Z_{s}^{i}dW_{s}.
\end{align*}
To conclude, $(X,Y,Z)$ is a solution for the FBSDE \eqref{eq:FBSDE_main} and use the same method in (9) and (10) we have $(Y,Z\cdot W)\in\mathcal{S}^{\infty}(\mathbb{R}^{n})\times BMO$.\\
\\
$\mathbf{Step\ 3}$: We show the minimality of the solution that we construct above.\\
Let us assume that $(\widetilde{X},\widetilde{Y},\widetilde{Z})$ is another solution of \eqref{eq:FBSDE_main} such that $(\widetilde{X},\widetilde{Y},\widetilde{Z}\cdot W)\in \mathcal{S}^{p}(\mathbb{R}^{n})\times \mathcal{S}^{\infty}(\mathbb{R}^{n}) \times BMO$ for any $p\geq2$.\\
From the assumption $(\mathscr{A}4)$, we have
\begin{align*}
-K(1+|y|+|z^{i}|^{2})\leq g^{i}(t,\widetilde{X},y,z),~\forall y\in\mathbb{R}^{n},~z^{i}\in\mathbb{R}^{d},~i=1,\cdots,n,~t\in[0,T]
\end{align*}
and we may conclude immediately that
\begin{align*}
Y_{t}^{(0)}\leq \widetilde{Y}_{t},~t\in[0,T]
\end{align*}
from \cite[Theorem 2.2]{19}. On the other hand, since $b$ is nondecreasing in $y$, from our construction we have
\begin{align*}
b(t,x,Y_{t}^{(0)})\leq b(t,x,\widetilde{Y}_{t}),~t\in[0,T].
\end{align*}
So by \cite[Theorem 1.1]{18} we have
\begin{align*}
X_{t}^{(0)}\leq \widetilde{X}_{t},~t\in[0,T]
\end{align*}
and then
\begin{align*}
&h^{i}(X_{T}^{(0)})\leq h^{i}(\widetilde{X}_{T}),\\
&g^{i}(t,X_{t}^{(0)},y,z)\leq g^{i}(t,\widetilde{X}_{t},y,z),~i=1,\cdots,n,~t\in[0,T].
\end{align*}
Again by \cite[Theorem 2.2]{19}, we conclude that $Y_{t}^{(1)}\leq\widetilde{Y}_{t},~t\in[0,T]$.\\
Iterating the procedure we get
\begin{align*}
&Y_{t}^{(0)}\leq Y_{t}^{(1)}\leq\cdots Y_{t}^{(k)}\leq\cdots\leq\widetilde{Y}_{t},\\
&X_{t}^{(0)}\leq X_{t}^{(1)}\leq\cdots X_{t}^{(k)}\leq\cdots\leq\widetilde{X}_{t},\\
&k\in\mathbb{N},~t\in[0,T].
\end{align*}
Therefore, $Y\leq\widetilde{Y},~X\leq\widetilde{X}$, that is the minimality of the solution $(X,Y)$.
\end{proof}

\begin{theorem}
  Assume $(x_{0},b,\sigma,h,g)$ and $(\overline{x}_{0},\overline{b},\sigma,\overline{h},\overline{g})$ satisfy $(\mathscr{A}1)-(\mathscr{A}7)$. 
  Let $(X,Y,Z)$ (resp. $(\overline{X},\overline{Y},\overline{Z})$) be the minimal solution of FBSDE \eqref{eq:FBSDE_main} associated to $(x_{0},b,\sigma,h,g)$ (resp. $(\overline{x}_{0},\overline{b},\sigma,\overline{h},\overline{g})$). If it holds that $x_{0}\leq \overline{x}_{0}$ and for any $i=1,\cdots,n,~t\in[0,T]$,
\begin{align*}
b^{i}(t,x,y)\leq \overline{b}^{i}(t,\overline{x},\overline{y})
\end{align*}
for any $x,\overline{x},y,\overline{y}\in\mathbb{R}^{n}$ satifying $x^{i}=\overline{x}^{i},x^{j}\leq\overline{x}^{j},j\neq i,y\leq\overline{y}$ and
\begin{align*}
&h(x)\leq \overline{h}(\overline{x})\\
&g^{i}(t,x,y,z^{i})\leq \overline{g}^{i}(t,\overline{x},\overline{y},z^{i})
\end{align*}
for any $z^{i}\in\mathbb{R}$ and $x,\overline{x},y,\overline{y}\in\mathbb{R}^{n}$ satifying $y^{i}=\overline{y}^{i},y^{j}\leq\overline{y}^{j},j\neq,x\leq\overline{x}$, we have
\begin{align*}
X_{t}\leq \overline{X}_{t},~Y_{t}\leq \overline{Y}_{t},~t\in[0,T].
\end{align*}
\end{theorem}
\begin{proof}
Since the solution constructed in Theorem 2.1 is the minimal one, we consider $(Y^{(0)},Z^{(0)})=(\overline{Y}^{(0)},\overline{Z}^{(0)})$ satisfying
\begin{align*}
Y_{t}^{i,(0)}=-C-C\int_{t}^{T}(1+|Y_{s}^{(0)}|+|Z_{s}^{(0)}|^{2})ds-\int_{t}^{T}Z_{s}^{i,(0)}dW_{s},~i=1,\cdots,n,~t\in[0,T].
\end{align*}
Then we consider
\begin{align*}
&X_{t}^{(0)}=x_{0}+\int_{0}^{t}b(s,X_{s}^{(0)},Y_{s}^{(0)})ds+\int_{0}^{t}\sigma(s,X_{s}^{(0)})dW_{s},\\
&\overline{X}_{t}^{(0)}=\overline{x}_{0}+\int_{0}^{t}\overline{b}(s,\overline{X}_{s}^{(0)},\overline{Y}_{s}^{(0)})ds+\int_{0}^{t}\sigma(s,\overline{X}_{s}^{(0)})dW_{s}.
\end{align*}
By our assumptions, it holds that
\begin{align*}
x_{0}^{i}\leq\overline{x}_{0}^{i},\quad b^{i}(t,x,Y_{t}^{(0)})\leq\overline{b}^{i}(t,\overline{x},\overline{Y}_{t}^{(0)}),~t\in[0,T]
\end{align*}
for any $x,\overline{x}\in\mathbb{R}^{n}$ satifying $x^{i}=\overline{x}^{i},x\leq\overline{x}^{j},j\neq i$, which implies
\begin{align*}
X_{t}^{(0)}\leq\overline{X}_{t}^{(0)},~t\in[0,T]
\end{align*}
due to the \cite[Theorem 1.1]{18}. Now we introduce $(Y^{(1)},Z^{(1)})$ and $(\overline{Y}^{(1)},\overline{Z}^{(1)})$ as follows
\begin{align*}
&Y_{t}^{i,(1)}=h^{i}(X_{T}^{(0)})+\int_{t}^{T}g^{i}(s,X_{s}^{(0)},Y_{s}^{(1)},Z_{s}^{i,(1)})ds-\int_{t}^{T}Z_{s}^{i,(1)}dW_{s},\\
&\overline{Y}_{t}^{i,(1)}=\overline{h}^{i}(\overline{X}_{T}^{(0)})+\int_{t}^{T}\overline{g}^{i}(s,\overline{X}_{s}^{(0)},\overline{Y}_{s}^{(1)},\overline{Z}_{s}^{i,(1)})ds-\int_{t}^{T}\overline{Z}_{s}^{i,(1)}dW_{s}.
\end{align*}
We can deduce
\begin{align*}
Y_{t}^{(1)}\leq \overline{Y}_{t}^{(1)},~t\in[0,T]
\end{align*}
from
\begin{align*}
h^{i}(X_{T}^{(0)})\leq\overline{h}^{i}(\overline{X}_{T}^{(0)}),~g^{i}(t,X_{t}^{(0)},y,z^{i})\leq\overline{g}^{i}(t,\overline{X}_{t}^{(0)},\overline{y},z^{i}),~i=1,\cdots,n,~t\in[0,T]
\end{align*}
for any $z^{i}\in\mathbb{R},y,\overline{y}\in\mathbb{R}^{n}$ satifying $y^{i}=\overline{y}^{i},y^{j}\leq\overline{y}^{j},j\neq i$ and \cite[Theorem 2.2]{19}.\\
We can claim that after iterating several times in the same way,
\begin{align*}
X_{t}^{(k)}\leq\overline{X}_{t}^{(k)},~Y_{t}^{(k)}\leq\overline{Y}_{t}^{(k)}
\end{align*}
holds for any $k\geq0,t\in[0,T]$.\\
Therefore we have
\begin{align*}
&X_{t}=\lim_{k\rightarrow\infty}X_{t}^{(k)}\leq\lim_{k\rightarrow\infty}\overline{X}_{t}^{(k)}=\overline{X}_{t},\\
&Y_{t}=\lim_{k\rightarrow\infty}Y_{t}^{(k)}\leq\lim_{k\rightarrow\infty}\overline{Y}_{t}^{(k)}=\overline{Y}_{t},~t\in[0,T].
\end{align*}
\end{proof}

\end{document}